\documentclass[12pt]{article}
\usepackage{amsmath}
\usepackage{amssymb}
\usepackage{graphicx}
\usepackage[cp1251]{inputenc}
\usepackage[russian]{babel}
\newcommand{\il}[2]{\int\limits_{#1}^{#2}}
\newcommand{\ilp}[1]{\int\limits_{#1}^{+\infty}}

\newcommand{\ph}{\phantom{a}}
\newcommand{\phh}{\phantom{aaa}}

\newcommand{\sist}[2]{\left\{
\begin{array}{l}
{#1}\\
\ph\\
{#2}
\end{array}
\right.}

\begin{document}

\vskip 20pt

MSC 34C10

\vskip 20pt

\centerline{\bf Oscillation and interval oscillation criteria for }
 \centerline{\bf linear matrix Hamiltonian systems}

\vskip 20 pt

\centerline{\bf G. A. Grigorian}
\centerline{\it Institute  of Mathematics NAS of Armenia}
\centerline{\it E -mail: mathphys2@instmath.sci.am}
\vskip 20 pt

\noindent
Abstract. We use the Riccati equation method with other ones to establish new oscillation and interval oscillation criteria for linear matrix Hamiltonian systems. We investigate the oscillation problem for linear matrix Hamiltonian systems in a new direction, which is to break the positive definiteness condition, imposed on one of the coefficients of the system.

\vskip 20 pt

Key words: Riccati equation, oscillation, interval oscillation, conjoined (prepared, preferred) solutions, unitary transformation, comparison theorem.

\vskip 20 pt

{\bf 1. Introduction.} Let $A(t)\ph B(t)$ and $C(t)$, be complex valued continuous matrix functions on $[t_0,+\infty)$ and let $B(t)$ and $C(t)$ be Hermitian, i.e., $B(t) = B^*(t), \ph C(t) =~ C^*(t), \linebreak t\ge t_0$ (here and after $*$ denotes the conjugation sign). Consider the linear matrix  Hamiltonian system
$$
\sist{\Phi'= A(t)\Phi + B(t)\Psi,}{\Psi' = C(t)\Phi - A^*(t)\Psi, \phh t\ge t_0,} \eqno (1.1)
$$
By a solution of this system we mean an ordered pair $(\Phi(t), \Psi(t))$ of  continuously differentiable matrix functions $\Phi(t)$ and $\Psi(t)$ of dimension $n\times n$ on $[t_0, +\infty)$, satisfying (1.1)  on $[t_0,+\infty)$.

{\bf Definition 1.1}. {\it A solution  $(\Phi(t), \Psi(t))$ of the system (1.1) is called conjoined (or prepared, preferred) if  $\Phi^*(t)\Psi(t) = \Psi^*(t)\Phi(t), \ph t\ge t_0$.}

{\bf Definition 1.2.} {\it A conjoined solution  $(\Phi(t), \Psi(t))$ of the system (1.1) is called oscilla- \linebreak tory if   $\det \Phi(t)$  has arbitrary large zeroes.}

{\bf Definition 1.3} {\it The system (1.1) is called oscillatory if its all conjoined solutions  are oscillatory}.

Let $[a,b]\subset [t_0,+\infty)$.

{\bf Definition 1.4.} {\it A conjoined solution  $(\Phi(t), \Psi(t))$ of the system (1.1) is called oscilla- \linebreak tory on the interval $[a,b]$  if   $\det \Phi(t)$  vanishes on $[a,b]$.}

{\bf Definition 1.5} {\it The system (1.1) is called oscillatory  on the interval $[a,b]$, if its all conjoined solutions  are oscillatory on $[a,b]$}.

Study of the oscillatory  behavior  of the system (1.1) is an important problem of qualitative theory of differential equations and many works are devoted to it (see \linebreak e.g., [1-12] and cited works therein). Usually the oscillation behavior of the system (1.1) is studied under the hypothesis that the matrix function $B(t)$ is positive definite  on $[t_0,+\infty)$, and this restriction  is essential from the point of view of the using methods of investigations. Meanwhile in the applications "the nature" \hskip 2pt of the restriction on $B(t)$ is that it must be non negative definite (The Legendre's condition).

 In [11] two oscillation criteria are obtained in a new direction which is to break the positive definiteness restriction imposed on $B(t)$. In [11] the last restriction was replaced by the non negative definiteness condition with the condition of solvability of the linear matrix equation
$$
\sqrt{B(t)} X [A(t) \sqrt{B(t)} - \sqrt{B(t)}'] = A(t) \sqrt{B(t)} - \sqrt{B(t)}', \phh t\ge t_0. \eqno (1.2)
$$

{\bf Remark 1.1.} {\it Eq. (1.2) has always a solution on $[t_0,+\infty)$ when $B(t)$ is invertible, in particular, when $B(t)$ positive definite on $[t_0,+\infty)$ ($X=\sqrt{B(t)}^{-1} \ph t\ge t_0$). But  it can also have a solution on $[t_0,+\infty)$ in some cases when $B(t)$ is not positive definite but it is nonnegative definite  (see [11]).}

Another replacements  of the mentioned above restriction are considered in [12], in which some new oscillation and interval oscillation criteria for the system (1.1) are obtained.

In this paper we continue the study of the oscillation problem of the system (1.1) in the mentioned above  direction. The Riccati equation method used to obtain new oscillation and interval oscillation criteria.  The unitary transformation approach allows  to obtain  oscillation and interval oscillation  criteria without solvability condition, imposed on Eq. (1.2).


{\bf 2. Main results}. The non negative (positive) definiteness of any Hermitian matrix we denote by $H\ge 0 (H > 0)$. Hereafter we will always assume that $B(t) \ge 0, \ph t\ge t_0$ (then $\sqrt{B(t)}, \ph t\ge t_0$ exists) and, when it is necessary, we will asumme that$\sqrt{B(t)}$ is  continuously differentiable on $[t_0,+\infty)$ (or an interval $[a,b] \subset [t_0,+\infty)$).

Let $F(t)$ be a matrix function of dimension $n\times n$ on $[t_0, +\infty)$. Set:
$$
A_F(t) \equiv F(t)[A(t)\sqrt{B(t)} - \sqrt{B(t)}'] = (a_{Fjk}(t))_{j,k =1}^n,
$$
$$
C_B(t) \equiv \sqrt{B(t)} C(t) \sqrt{B(t)} =(c_{Bjk}(t))_{j,k =1}^n,
$$
$$
\theta_{Fj}(t) \equiv c_{Bjj}(t) + \sum\limits_{\stackrel{m =1} {m\ne j}}^n |a_{Fmj}(t)|^2, \ph j= \overline{1,n}, \ph t \ge t_0.
$$

{\bf Theorem 2.1.} {\it Let the following conditions be satisfied.

1) Eq. (1.2) has a solution $F(t)$ on $[t_0,+\infty)$;

2)  for some $j \in \{1, ..., n\}$  the scalar equation
$$
\phi'' + 2 \mathfrak{Re}\hskip 2pt a_{Fjj}(t) \phi' + \theta_{Fj}(t) \phi = 0, \phh t\ge t_0 \eqno (2.1)
$$

is oscillatory.

Then the system (1.1) is also oscillatory.}

$\phantom{aaaaaaaaaaaaaaaaaaaaaaaaaaaaaaaaaaaaaaaaaaaaaaaaaaaaaaaaaaaaaaaaaaaa} \Box$

{\bf Theorem 2.2.} {\it Let the following conditions be satisfied.

1') Eq. (1.2) has a solution $F(t)$ on $[a,b]$;

2')  for some $j \in \{1, ..., n\}$  the scalar equation
$$
\phi'' + 2 \mathfrak{Re}\hskip 2pt a_{Fjj}(t) \phi' + \theta_{Fj}(t) \phi = 0, \phh t\in [a,b]
$$

is oscillatory on $[a,b]$.

Then the system (1.1) is also oscillatory on $[a,b]$.}

$\phantom{aaaaaaaaaaaaaaaaaaaaaaaaaaaaaaaaaaaaaaaaaaaaaaaaaaaaaaaaaaaaaaaaaaaa} \Box$

{\bf Remark 2.1.} {\it An explicit interval oscillation criterion for second order linear ordinary differential equations (therefore for Eq. (2.1)) id obtained in [13] (see [13], Theorem 3.2)}

\vskip 10pt

The next result is based on the use of an unitary transformation, which allows us to overcome the restriction of solvability of Eq. (1.2), presented in the conditions of \linebreak Theorem~ 2.1.

Let $p_{jk}(t), \ph j,k =1,2$ be real-valued locally integrable functions on $[t_0,+\infty)$. Consider the linear system of ordinary differential equations
$$
\sist{\phi' = p_{11}(t) \phi + p_{12}(t) \psi,}{\psi' = p_{21}(t) \phi + p_{22}(t) \psi, \ph t\ge t_0.} \eqno (2.1)
$$

{\bf Definition 2.1.} {\it A solution $(\phi(t), \psi(t))$ of the system (2.1) is called oscillatory if $\phi(t)$ has arbitrary large zeroes.}

{\bf Definition 2.2.} {\it The system (2.1) is called oscillatory if its all solutions are oscillatory.}

{\bf Definition 2.3.} {\it  A solution $(\phi(t), \psi(t))$ of the system (2.1) is called oscillatory  on the interval $[a,b]$, if $\phi(t)$ vanishes on $[a,b]$.}

{\bf Definition 2.4.} {\it The system (2.1) is called oscillatory on the interval $[a,b]$  if its all solutions are oscillatory on  $[a,b]$.}

Let $U_B(t)$ be an unitary matrix function of dimension $n\times n$ on $[t_0,+\infty)$ such that
$$
B(t) = U_B^*(t) B_0(t) U_B(t), \phh t\ge t_0., \eqno (2.2)
$$
where $B_0(t) \equiv diag\{b_1(t), ..., b_n(t)\}, \ph t\ge t_0$ - is a diagonal matrix function on $[t_0, +\infty)$.

{\bf Remark 2.2.} {\it It is well known that for any Hermitian matrix $H$ of dimension $n\times n$ there exists an unitary matrix (transformation) $U_H$ such that $H = U_H^* \hskip 2pt diag\{h_1, ..., h_n\} U_H,$ where $h_1, ..., h_n$ are real numbers.}

Hereafter we will assume that $U_B(t)$ is continuously differentiable on $[t_0,+\infty)$ and $B_0(t)$ is continuous on $[t_0, +\infty)$. Set:
$$
A_B^0(t) \equiv U_B(t) [A(t) U_B^*(t) - \{U_B^*(t)\}'] = (a_{jk}^0(t))_{jk=1}^n,
$$
$$
C_B^0(t) \equiv U_B(t) C(t) U_B^*(t) = (c_{jk}^0(t))_{jk=1}^n,
$$
$$
\biggl[\frac{|a_{mj}^0(t)|^2}{b_m(t)}\biggr]_0 \equiv \sist{\frac{|a_{mj}^0(t)|^2}{b_m(t)},\ph \mbox{if} \ph b_m(t) \ne 0,}{0, \phh  \mbox{if} \phh b_m(t) = 0,} \phh m= \overline{1,n},
$$
$$
\chi_j(t) \equiv c_{jj}^0(t) + \sum\limits_{\stackrel{m =1} {m\ne j}}^n\biggl[\frac{|a_{mj}^0(t)|^2}{b_m(t)}\biggr]_0,\phh j=\overline{1,n}, \phh t\ge t_0.
$$

\vskip 20pt

{\bf Theorem 2.3.} {\it Let the following conditions be satisfied:

3) $b_m(t) \ge 0, \ph m= \overline{1,n}, \ph t\ge t_0;$

4) for some $j \in \{1,...,n\}$  the function $ \chi_j(t)$ is continuous on $[t_0,+\infty)$ and  the scalar system
$$
\sist{\phi' = 2 \hskip 2pt \mathfrak{Re}\hskip 2pt a_{jj}^0(t) \phi + b_j(t) \psi,}{\psi' =  \chi_j(t) \phi, \phh t\ge t_0} \eqno (2.3)
$$
is oscillatory.

Then the system (1.1) is also oscillatory.}

$\phantom{aaaaaaaaaaaaaaaaaaaaaaaaaaaaaaaaaaaaaaaaaaaaaaaaaaaaasssssssssssssaaaaa} \Box$

\vskip 10pt

{\bf Remark 2.3}. {\it Explicit oscillatory criteria for the system (2.1) (therefore for the system (2.3)) are obtained in [14].}

\vskip 10pt

{\bf Corollary 2.1.} {\it Let the following conditions be satisfied

\noindent
5) $B(t) = diag \hskip 3pt \{b_1(t), \dots, b_n(t)\}, \ph b_m(t) \ge 0, \ph m=\overline{1,n}, \ph t \ge t_0$,

\noindent
6) for for some $j\in \{1,\dots,n\}$
$$
\ilp{t_0}b_j(\tau) d \tau = -\ilp{t_0}c_{jj}(\tau) d \tau = +\infty.
$$
Then the system
$$
\sist{\Phi' = \phh B(t) \Psi,}{\Psi' = C(t) \Phi, \phh t \ge t_0} \eqno (2.4)
$$
is oscillatory.}

\vskip 10pt

{\bf Remark 2.4.} {\it Corollary 2.1 is a generalization of Leighton's oscillation criterion (see [15, Theorem 2.24]).}

\vskip 10pt

{\bf Theorem 2.4.} {\it Let the following conditions be satisfied:

3') $b_m(t) \ge 0, \ph m= \overline{1,n}, \ph t\in [a,b];$

4') for some $j \in \{1,...,n\}$  the function $ \chi_j(t)$ is continuous on $[a,b]$ and  the scalar system
$$
\sist{\phi' = 2 \hskip 2pt \mathfrak{Re}\hskip 2pt a_{jj}^0(t) \phi + b_j(t) \psi,}{\psi' =  \chi_j(t) \phi, \phh t\in [a,b]}
$$
is oscillatory on $[a,b]$.

Then the system (1.1) is also oscillatory on $[a,b]$.}

$\phantom{aaaaaaaaaaaaaaaaaaaaaaaaaaaaaaaaaaaaaaaaaaaaaaaaaaaaaaaaaaaaaaaaaaaaaa} \Box$

\vskip 10pt

{\bf Corollary 2.2.} {\it Let the following conditions be satisfied

5') $B(t) = diag \hskip 3pt \{b_1(t), \dots, b_n(t)\}, \ph b_m(t) \ge 0, \ph m=\overline{1,n}, \ph t\in [a,b]$,

6') for for some $j\in \{1,\dots,n\}$
$$
\int\limits_a^b \min\bigl[b_j(t)
, - c_{jj}(t)\bigr]d t
 \ge \pi.
$$

Then the system (2.4) is oscillatory on the interval $[a,b]$}

\vskip 10pt

{\bf 3. Proof of the main results}. Let $f_k(t), \ph g_k(t)$ and $h_k(t), \ph k=1,2$ be real-valued  continuous functions on $[t_0,+\infty)$. Consider the scalar Riccati equations
$$
y' + f_k(t) y^2 + g_k(t) y + h_k(t) = 0, \ph t\ge t_0, \ph k=1,2 \eqno (3.1_k)
$$
and the differential inequalities
$$
\eta' + f_k(t) \eta^2 + g_k(t) \eta + h_k(t) = 0, \ph t\ge t_0, \ph  k=1,2. \eqno (3.2_k)
$$

{\bf Remark 3.1.} {\it Every solution of Eq. $(3.1_k)$ on $[t_1,t_2) \ph (t_0 \le t_1< t_2 \le +\infty)$ is also a solution of the inequality $(3.2_k), \ph k =1,2.$}

{\bf Remark 3.2.} {\it If $f_k(t) \ge 0, \ph t\ge t_0,$ then every solution of the linear equation
$$
\zeta' + g_k(t) \zeta + h_k(t) = 0, \phh t\ge t_0
$$
is also a solution of the inequality $(3.2_k), \ph k =1,2.$}

The following comparison theorem plays a crucial role in thee proof of the main results.

{\bf Theorem 3.1 [16, Theorem 3.1].} {\it Let Eq. $(3.1_2)$ have a real valued solution $y_2(t)$ on $[t_0,\tau_0) \ph (t_0 < \tau_0 \le +\infty)$ and let the following conditions be satisfied: $f_1(t) \ge 0$  and $\il{t_0}{t} \exp\biggl\{\il{t_0}{\tau}[f(s)(\eta_1(s) + \eta_2(s)) + g(s)]d s\biggr\}[(f_1(\tau) - f(\tau))y_2^2(\tau) + (g_1(\tau) - g(\tau)) y_2(\tau) + h_1(\tau) - h(\tau)]d \tau \ge 0, \ph t\in [t_0,\tau_0)$ where $\eta_1(t)$ and $\eta_2(t)$ are solutions of the inequalities $(3.2_1)$ and $(3.2_2)$ respectively on $[t_0, \tau_0)$ such that $\eta_j(t_0) \ge y_2(t_0), \ph j=1,2.$ Then for every $\gamma_0 \ge y_2(t_0)$ Eq. $(3.1_1)$ has a solution $y_1(t)$ on $[t_0, \tau_0)$, satisfying the condition $y_1(t_0) = \gamma_0.$}

{\bf Remark 3.3.} {\it One can easily verify, that in the case $\tau_0 <+\infty$ Theorem 3.1 remains valid if we replace $[t_0,\tau_0)$ by $[t_0,\tau_0]$ in it.}

Set $E(t) \equiv p_{11}(t) - p_{22}(t), \ph t \ge t_0$.

{\bf Theorem 3.2 [14, Theorem 2.4].} {\it Let the following conditions be satisfied:

\noindent
 $p_{12}(t) \ge 0, \ph t\ge t_0$;

\noindent
$\int\limits_{t_0}^{+\infty} p_{12}(t)\exp\bigl\{-\int\limits_{t_0}^tE(\tau) d\tau\bigr\}= - \int\limits_{t_0}^{+\infty}p_{21}(t)\exp\bigl\{\int\limits_{t_0}^tE(\tau) d\tau\bigr\} d t = +\infty$.

\noindent
Then the system (2.1) is oscillatory.} $\Box$

{\bf Theorem 3.3 [14, Theorem 2.3].} {\it Let the following conditions be satisfied:

\noindent
 $p_{12}(t) \ge 0, \ph t\in [a;b];$

\noindent
 $\int\limits_a^b \min\biggl[p_{12}(t)\exp\bigl\{-\int\limits_a^tE(\tau) d\tau\bigr\}, - p_{21}(t)\exp\bigl\{\int\limits_a^tE(\tau) d\tau\bigr\}\biggr]d t \ge \pi.$

\noindent
Then the system (2.1) is oscillatory on $[a;b]$.} $\Box$

Consider the scalar Riccati equation
$$
y' + p_{12}(t) y^2 + (p_{11}(t) - p_{22}(t)) y - p_{21}(t) = 0, \phh t\ge t_0.
$$
The solutions $y(t)$  of this equation, existing on some interval $[t_1,t_2) (t_0 \le t_1 < t_2 \le + \infty)$ are connected with solutions $(\phi(t), \psi(t))$ of the system (2.1) by relations (see [16])
$$
\phi(t) = \phi(t_1) \exp\biggl\{\il{t_1}{t}[p_{12}(\tau) y(\tau) + a_{11}(\tau)]d\tau\biggr\}, \ph \phi_1(t_1) \ne 0, \ph \psi(t) = y(t) \phi(t),  \eqno (3.3)
$$
$t\in [t_1.t_2)$.

Let $p(t)$ and $q(t)$ be real-valued locally integrable functions on $[t_0,+\infty)$. Consider the second order linear ordinary differential equation
$$
\phi'' + p(t)\phi' + q(t) \phi = 0, \phh t\ge t_0 \eqno (3.4)
$$
and the corresponding scalar Riccati one
$$
y' + y^2 + p(t) y + q(t) = 0, \phh t\ge t_0. \eqno (3.5)
$$
Since Eq. (3.4) is equivalent to the system
$$
\sist{\phi' = \ph \psi,}{\psi' = - q(t) \phi - p(t) \psi, \ph t\ge  t_0}
$$
by (3.3) we have that the solutions $y(t)$ of Eq. (3.5), existing on an interval $[t_1,t_2)$, are connected with solutions $\phi(t)$ of Eq. (3.4) by relations
$$
\phi(t) = \phi(t_1)\exp\biggl\{\il{t_1}{t}\Bigl[y(\tau) + p(\tau)\Bigr] d \tau\biggr\}, \phh \phi(t_1) \ne 0, \phh t\in [t_1,t_2). \eqno (3.6)
$$

Consider the matrix Riccati equation
$$
Z' + Z B(t) Z + A^*(t) Z + Z A(t) - C(t) = 0, \phh t\ge t_0. \eqno (3.7)
$$
It is not difficult to verify that the solutions $Z(t)$ of this equation, existing on an interval $[t_1, t_2) \ph (t_0 \le t_1 < t_2 \le + \infty)$ are connected with solutions $(\Phi(t), \Psi(t))$ of the system (1.1) by the relations
$$
\Phi'(t) = [A(t) + B(t) Z(t)] \Phi(t), \ph \Phi(t_1) \ne 0, \ph \Psi(t) = Z(t) \Phi(t), \ph t\in [t_1,t_2). \eqno (3.8)
$$

{\bf 3.1. Proof of Theorem 2.1.} Suppose the system (1.1) is not oscillatory. Then there exists a conjoined solution
$(\Phi(t), \Psi(t))$ of that system such that $det \hskip 2pt \Phi(t) \ne 0, \ph t\ge t_1$ for some $t_1 \ge t_0$. By (3.8)(3.9) from here it follows that $Z(t) \equiv \Psi(t)\Phi^{-1}(t), \ph t\ge t_1$ is a Hermitian solution  of Eq. (3.7)(3.8) on $[t_1,+\infty),$ i. e., $Z^*(t) = Z(t)$ and
$$
Z'(t) + Z(t) B(t) Z(t) + A^*(t) Z(t) + Z(t) A(t) - C(t) = 0, \phh t\ge t_1.
$$
Multiply both sides of this equality at left and at right by $\sqrt{B(t)}, \ph t\ge t_1.$ Taking into account the equality
$$
\sqrt{B(t)} Z'(t) \sqrt{B(t)} = [\sqrt{B(t)} Z(t) \sqrt{B(t)}]' - \sqrt{B(t)}' Z(t) \sqrt{B(t)} - \sqrt{B(t)} Z(t) \sqrt{B(t)}', \ph t\ge t_1
$$
and the condition 1) of the theorem we obtain
$$
V'(t) + V^2(t) + A_F^*(t) V(t) + V(t) A_F(t) - C_B(t) = 0, \phh t\ge t_1, \eqno (3.9)
$$
where $V(t) \equiv \sqrt{B(t)} Z(t) \sqrt{B(t)}, \ph t\ge t_1$. Denote by $[M]_{jk}$ the $j k$ -th entry of any square matrix $M \ph (j, k =\overline{1,n})$. Set: $[V(t)]_{jk} \equiv v_{jk}(t), \ph t\ge t_1, \ph k=\overline{1,n}.$
Since $V(t)$ is a Hermitian matrix function on $[t_1, +\infty)$ it is not difficult to verify that
$$
[V^2(t)]_{11} = v_{11}^2(t) + |v_{12}(t)|^2 + ... + |v_{1 n}(t)|^2,
$$
$$
[V^2(t)]_{22} = |v_{21}(t)|^2  +  v_{22}^2(t) +  ... + |v_{2 n}(t)|^2,
$$
$$
 - - - - - - - - - - - - - - - - - - -
 $$
$$
[V^2(t)]_{nn} = |v_{n1}(t)|^2 + |v_{n 2}(t)|^2 + ... + v_{n n}^2(t)
$$
$$
[V(t) A_F(t)]_{j j} = \sum\limits_{m=1}^n v_{jm}(t) a_{Fmj}(t), \ph [A_F^*(t) V(t)]_{j j} = \sum \limits_{m=1}^n\overline{v_{jm}(t)} \hskip 4pt\overline{a_{Fmj}(t)}, \ph t\ge t_1.
$$
From here and from the equalities $v_{j m}(t) = \overline{v_{m j}(t)}, \ph m =\overline{1,n}, \ph t\ge t_1$ we obtain
$$
v_{j j}'(t) + v_{j j}^2(t) + 2 \mathfrak{Re} \hskip 3pt a_{Fj j}(t)  v_{j j}(t) +  \sum\limits_{\stackrel{m=1} {m \ne j}}^n |v_{j m}(t) + \overline{a_{Fm j}(t)}|^2 - \theta_{Fj}(t) = 0, \phh t\ge t_1. \eqno (3.10)
$$
Consider the scalar Riccati equations
$$
y' + y^2 + 2 \mathfrak{Re} \hskip 3pt a_{F j j}(t) y - \theta_{F j}(t) = 0, \phh t\ge t_1, \eqno (3.11)
$$

$$
y' + y^2 + 2 \mathfrak{Re} \hskip 3pt a_{Fj j}(t) y - \theta_{Fj}(t) + \sum\limits_{\stackrel{m=1} {m \ne j}}^n |v_{j m}(t) + \overline{a_{Fm j}(t)}|^2  = 0, \phh t\ge t_1. \eqno (3.12)
$$
By (3.10) $v_{j j}(t)$ is a solution to the last equation on $[t_1, + \infty)$. Since \linebreak $\sum\limits_{\stackrel{m=1} {m \ne j}}^n |v_{j m}(t) + \overline{a_{Fm j}(t)}|^2 \ge~ 0, \ph t\ge t_1$, using Theorem 3.1 to the pair of the equations (3.11) and (3.12) we conclude that Eq. (3.11) has a solution $y_1(t)$ on $[t_1, + \infty)$. Then by (3.6) $\phi_1(t) \equiv  \exp\biggl\{\il{t_1}{t}\Bigl[y_1(\tau) + \mathfrak{Re} \hskip 3pt a_{Fj j}(t) y_1(\tau)\Bigr] d \tau\biggr\}, \ph t \ge t_1$ is a solution of Eq. (2.1) on $[t_1,+ \infty)$, which can be continued on $[t_0, +\infty)$ as a solution of Eq. (2.1). Since $\phi_1(t) > 0, \ph t \ge t_1$ Eq. (2.1) is not oscillatory, which contradicts the condition 2) of the theorem. The obtained contradiction completes the proof of the theorem.

{\bf Remark 3.4.} {\it  Theorem 2.2 can be proved by analogy of the proof of Theorem 2.1 by taking into account Remark 3.3.}

{\bf 3.2. Proof of Theorem 2.3}. Suppose the system (1.1) is not oscillatory. Then there exists a conjoined solution $(\Phi(t), \Psi(t))$ of (1.1) such that $det \hskip 3pt \Phi(t) \ne~ 0, \ph t\ge t_1$ for some $t_1 \ge t_0$. By virtue of (3.5) from here it follows that $Z(t) \equiv \Psi(t) \Phi^{-1}(t), \ph t\ge t_1$  is a Hermitian solution of Eq. (3.7)(3.4) on $[t_1, +\infty)$, that is  $Z^*(t) = Z(t), \ph t\ge t_1$ and
$$
Z'(t) + Z(t) B(t) Z(t) + A^*(t) Z(t) + Z(t) A(t) - C(t) = 0, \phh t\ge t_1
$$
Multiply both sides of the last equality at left by $U_B(t)$ and at right by $U_B^*(t)$. Taking into account (2.2) and the equality
$$
U_B(t)Z'(t) U_B^*(t) = [U_B(t)Z(t) U_B^*(t)]' - U'_B(t)Z(t) U_B^*(t) - U_B(t)Z'(t)[ U_B^*(t)]', \phh t\ge t_1,
$$
we obtain
$$
V'(t) + V(t) B_0(t) V(t) + [A^0_B(t)]^* V(t) + V(t) A_B^0(t) - C_B^0(t) = 0, \phh t\ge t_1, \eqno (3.13)
$$
where $V(t) \equiv U_B(t) Z(t) U_B^*(t), \ph t\ge t_1$.  Let $V(t) \equiv (v_{jk}(t))_{j,k =1}^n, \ph t\ge t_1$. Since $V(t)$ is a Hermitian matrix function it is not difficult to verify that
$$
[V(t) B_0(t) V(t)]_{11} = b_1(t) v_{11}^2(t) + b_2(t) |v_{12}(t)|^2 + ... + b_n(t)|v_{1n}(t)|^2,
$$
$$
[V(t) B_0(t) V(t)]_{22} = b_1(t) |v_{21}(t)|^2  +  b_2(t) v_{22}^2(t)  + ... + b_n(t)|v_{2n}(t)|^2,
$$
$$
- -  - - - - - - - - - - - - - - - -  - - -   - - - - - -  -
$$
$$
[V(t) B_0(t) V(t)]_{nn} =  b_1(t) |v_{n1}(t)|^2 + b_2(t)|v_{n2}(t)|^2 + ... + b_n(t) v_{nn}^2(t),
$$
$$
[V(t) A_B^0(t)]_{jj} = \sum\limits_{m=1}^n v_{jm}(t) a_{mj}(t), \phh [(A_B^0(t))^* V(t)]_{jj} = \sum\limits_{m=1}^n\overline{v_{jm}(t)}\hskip 2pt \overline{a_{mj}^0(t)}, \phh t\ge t_1.
$$
Taking into account the equalities $v_{jm}(t) = \overline{v_{mj}(t)}, \ph m = \overline{1,n}, \ph t\ge t_1$ from here we obtain
$$
v_{jj}'(t) + b_j(t) v_{jj}^2(t) + 2 \mathfrak{Re} \hskip 2pt a_{jj}^0(t)v_{jj}(t) + \sum\limits_{\stackrel{m=1} {m \ne j}}^n b_m(t)\biggl|v_{jm}(t) + \frac{\overline{a_{mj}^0}(t)}{b_m(t)}\biggr|_0^2 - \chi_j(t) = 0,  \eqno (3.14)
$$
$t \ge t_1$, where
$$
\biggl|v_{jm}(t) + \frac{\overline{a_{mj}^0}(t)}{b_m(t)}\biggr|_0 \equiv \sist{\biggl|v_{jm}(t) + \frac{\overline{a_{mj}^0}(t)}{b_m(t)}\biggr|, \ph \mbox{if} \ph b_m(t) \ne 0,}{0, \phh \mbox{if} \phh b_m(t) = 0,} \phh m=\overline{1,n}, \ph t\ge t_1.
$$
Consider the scalar Riccati equations
$$
y' + b_j(t) y^2 + 2 \mathfrak{Re} \hskip 2pt a_{jj}^0(t) y - \chi_j(t) = 0, \phh t \ge t_1, \eqno (3.15)
$$
$$
y' + b_j(t) y^2 + 2 \mathfrak{Re} \hskip 2pt a_{jj}^0(t) y + \sum\limits_{\stackrel{m=1} {m \ne j}}^n b_m(t)\biggl|v_{jm}(t) + \frac{\overline{a_{mj}^0}(t)}{b_m(t)}\biggr|_0^2 - \chi_j(t) = 0, \ph t\ge t_1. \eqno (3.16)
$$
By (3.14) $v_{jj}(t)$ is a solution to the last equation on $[t_1, + \infty)$. From the condition 4) of the theorem it follows that
$
\sum\limits_{\stackrel{m=1} {m \ne j}}^n b_m(t)\biggl|v_{jm}(t) + \frac{\overline{a_{mj}^0}(t)}{b_m(t)}\biggr|_0^2 \ge 0, \phh t\ge t_1.
$
Then using Theorem~ 3.1 to the pair of equations (3.15) and (3.16) we conclude that Eq. (3.15) has a solution $y(t)$ on $[t_1, + \infty)$. Hence in virtue of (3.1) the functions
$$
\phi(t) \equiv \exp\biggl\{\il{t_1}{t}[b_j(\tau) y(\tau) + 2 \mathfrak{Re} \hskip 2pt a_{jj}^0(\tau)]d\tau\biggr\}, \phh \psi(t) \equiv y(t) \phi(t), \phh t\ge t_1
$$
form a solution $(\phi(t), \psi(t))$ of the system (2.3) on $[t_1, + \infty)$, which can be continued on $[t_0, +\infty)$ as a solution of the system (2.3). Since, obviously, $\phi(t) > 0, \ph t\ge t_1$ the system (2.3) is not oscillatory, which contradicts the condition 4) of the theorem. The obtained contradiction completes the proof of the theorem.

{\bf Remark 3.5.} {\it  Theorem 2.4 can be proved by analogy of the proof of Theorem 2.3 by taking into account Remark 3.3.}

\vskip 10 pt

{\bf 3.3. Proof of Corollary 2.1.} Since according to the condition 5) $B(t)$ is a diagonal matrix, we can take the unitary transformation  $U_B(t) \equiv I$. Then  for the system (2.4) we will have $a_{jj}^0(t) \equiv 0, \ph \chi_j(t) = c_{jj}(t), \ph t \ge t_0.$ Then by Theorem 2.3 from the condition 5) it follows that the system (2.4) is oscillatory provided the scalar system
$$
\sist{\phi' = \phh b_j(t) \psi,}{\psi' = c_{jj}(t), \phh t \ge t_0}
$$
is oscillatory. By Theorem 3.2 this condition holds provided the condition 6) is satisfied. The corollary is proved.

Corollary 2.2 can be proved by analogy of the proof of Corollary 2.1 using Theorem~ 3.3 instead of Theorem 3.2.

\vskip 20 pt

\centerline{ \bf References}

\vskip 20pt

\noindent
1. L. Li, F. Meng and Z. Zheng, Oscillation results related to integral averaging technique\linebreak \phantom{a} for linear Hamiltonian systems, Dynamic Systems  Appli. 18 (2009), \ph \linebreak \phantom{a} pp. 725 - 736.

\noindent
2. F. Meng and  A. B. Mingarelli, Oscillation of linear Hamiltonian systems, Proc. Amer.\linebreak \phantom{a} Math. Soc. Vol. 131, Num. 3, 2002, pp. 897 - 904.

\noindent
3. Q. Yang, R. Mathsen and S. Zhu, Oscillation theorems for self-adjoint matrix \linebreak \phantom{a}   Hamiltonian
 systems. J. Diff. Equ., 19 (2003), pp. 306 - 329.

\noindent
4. Z. Zheng and S. Zhu, Hartman type oscillatory criteria for linear matrix Hamiltonian  \linebreak \phantom{a} systems. Dynamic  Systems  Appli., 17 (2008), pp. 85 - 96.

\noindent
5. Z. Zheng, Linear transformation and oscillation criteria for Hamiltonian systems. \linebreak \phantom{a} J. Math. Anal. Appl., 332 (2007) 236 - 245.

\noindent
6. I. S. Kumary and S. Umamaheswaram, Oscillation criteria for linear matrix \linebreak \phantom{a} Hamiltonian systems, J. Differential Equ., 165, 174 - 198 (2000).

\noindent
7. Sh. Chen, Z. Zheng, Oscillation criteria of Yan type for linear Hamiltonian systems, \linebreak \phantom{a} Comput.  Math. with Appli., 46 (2003), 855 - 862.

\noindent
8. Y. G. Sun, New oscillation criteria for linear matrix Hamiltonian systems. J. Math. \linebreak \phantom{a} Anal. Appl., 279 (2003) 651 - 658.

\noindent
9. K. I. Al - Dosary, H. Kh. Abdullah and D. Husein. Short note on oscillation of matrix \linebreak \phantom{a} hamiltonian systems. Yokohama Math. J., vol. 50, 2003.

\noindent
10. G. A. Grigorian, Oscillatory and Non Oscillatory criteria for the systems of two \linebreak \phantom{aa}   linear first order two by two dimensional matrix ordinary differential equations. \linebreak \phantom{aa}   Arch.  Math., Tomus 54 (2018), PP. 189 - 203.

\noindent
11. G. A. Grigorian, Oscillation criteria for linear matrix Hamiltonian systems. \linebreak \phantom{aa} Proc. Amer. Math. Sci, Vol. 148, Num. 8 ,2020, pp. 3407 - 3415.

\noindent
12. G. A. Grigorian.   Interval oscillation criteria for linear matrix Hamiltonian systems,\linebreak \phantom{a} vol. 50 (2020), No. 6, 2047–2057

\noindent
13.	G. A. Grigoryan, Some properties of solutions of second-order linear ordinary differential \linebreak \phantom{aa} equations, Trudy Inst. Mat. i Mekh. UrO RAN, 19:1 (2013),  69–80.

\pagebreak

\noindent
14.  G. A. Grigorian. Oscillatory criteria for the systems of two first - order Linear \linebreak \phantom{a} ordinary differential equations. Rocky Mount. J. Math., vol. 47, Num. 5,
 2017, \linebreak \phantom{a}  pp. 1497 - 1524

\noindent
15. C. A. Swanson. Comparison and oscillation theory of linear differential equations.   \linebreak \phantom{a} Academic press. New York and London, 1968.

\noindent
16. G. A. Grigorian,  On two comparison tests for second-order linear  ordinary\linebreak \phantom{aa} differential equations (Russian) Differ. Uravn. 47 (2011), no. 9, 1225 - 1240; trans-\linebreak \phantom{aa} lation in Differ. Equ. 47 (2011), no. 9 1237 - 1252, 34C10.


\end{document}